\documentclass[12pt]{amsart}
\usepackage{amsmath,amssymb,amsbsy,amsfonts,amsthm,latexsym,amsopn,amstext,
                                                       amsxtra,euscript,amscd}
\begin{document}

\title[generalized harmonic and oscillatory numbers]{On the properties of generalized \\ harmonic and oscillatory numbers. \\Simple proof of the Prime Number Theorem}

\author[R. M.~Abrarov]{R. M.~Abrarov}
\address{University of Toronto, Canada}
\email{rabrarov@physics.utoronto.ca}

\author[S. M.~Abrarov]{S. M.~Abrarov}
\address{York University, Toronto, Canada}
\email{abrarov@yorku.ca}

\date{\today}

\begin{abstract} 
We derived the sum identities for generalized harmonic and corresponding oscillatory numbers for which a sieve procedure can be applied. The obtained results enable us to understand better the properties of these numbers and their asymptotic behavior. On the basis of these identities a simple proof of the Prime Number Theorem is represented.
\end{abstract}

\maketitle

\section{Generalized harmonic numbers}

In our earlier report we have discussed the regular parts for basic functions of prime numbers \cite{RAbrarov07}. Before considering their oscillatory parts, we would like to discuss some important properties of the generalized harmonic and corresponding oscillatory numbers.

The generalized harmonic number in power \textit{s} is given by
\begin{equation} \label{Eq_1}
H_{x} \left(s\right)=\sum _{k=1}^{x}\frac{1}{k^{s} }  ,
\end{equation} 
where \textit{s} is any complex number. The basic properties of these numbers can be found elsewhere \cite{Wolfram}. Throughout this paper we use repeatedly M\"obius inversion formula \cite{Hardy79} and here we give it for references: 
if for all positive \textit{x} satisfied
\begin{equation} \label{Eq_2}
{\rm {\mathcal G}}\left(x\right) =\sum _{k=1}^{x}{\rm {\mathcal F}}\left(\frac{x}{k} \right)\\ 
\end{equation}
\text{then}
\begin{equation} \label{Eq_3}
{\rm {\mathcal F}}\left(x\right) =\sum _{k=1}^{x}\mu \left(k\right){\rm {\mathcal G}}\left(\frac{x}{k} \right)
\end{equation}
and vice versa, where
$$
\mu \left(n\right)=\mu _{n} =
\begin{cases} 1 & \text{if $n=1$,} 
\\ 
(-1)^{m} & \text{if $n$ is a product of $m$ distinct primes,}
\\ 
0 &\text{if the square of primes divides $n$.}
\end{cases}
$$
is M\"obius function.

From the definition \eqref{Eq_1} and M\"obius inversion formula \eqref{Eq_2} and \eqref{Eq_3} directly follows
$$
x^{s} \cdot H_{x} \left(s\right)=\sum _{k=1}^{x}\left(\frac{x}{k} \right)^{s}
$$
and
$$
x^{s} =\sum _{k=1}^{x}\mu _{k} \left(\frac{x}{k} \right) ^{s} H_{\frac{x}{k} } \left(s\right).
$$ 
Hence we get the important sum identity
\begin{equation}\label{Eq_4}
\sum_{k=1}^{x}\frac{\mu _{k} }{k^{s} }  H_{\frac{x}{k} } \left(s\right)=1.
\end{equation}
Using Stieltjes integration method, we can rewrite the same equation in integral form
\begin{equation}\label{Eq_5}
\int _{1-}^{x}H_{\frac{x}{y} } \left(s\right) \cdot dM_{y} \left(s\right)=1, 
\end{equation}
where $M_{y} \left(s\right)=\sum _{k=1}^{y}\frac{\mu _{k} }{k^{s} }  $ is corresponding oscillatory number in power \textit{s }(our notations are similar to those of commonly accepted \cite{Bateman04,Montgomery06,Tenenbaum95,Tenenbaum01} for the case \textit{s }= 0 and \textit{s }= 1, see below definitions \eqref{Eq_35}-\eqref{Eq_37}).

At \textit{s }=1 for ordinary harmonic number $H_{x} \left(1\right)\equiv H_{x} $, we have
\begin{equation}\label{Eq_6}
\sum _{k=1}^{x}\frac{\mu _{k} }{k}  H_{\frac{x}{k} } =1
\end{equation}
with corresponding integral form
\begin{equation}\label{Eq_7}
\int _{1-}^{x}H_{\frac{x}{y} }  \cdot dM_{y} =1.
\end{equation}

Applying in \eqref{Eq_6} the asymptotic formula for harmonic number
\[H_{\frac{x}{k} } =\sum _{n=1}^{x/k}\frac{1}{n}  =\log \left[\frac{x}{k}\right] +\gamma +\frac{1}{2\left[\frac{x}{k}\right]} -\frac{1}{12\left[\frac{x}{k}\right]^{2} } +\frac{1}{120\left[\frac{x}{k}\right]^{4} } -... ,\] 
where $\gamma = 0.5772156\: \dots $  is Euler's constant, and substituting $\left[\frac{x}{k}\right]$ approximatly by $\frac{x}{k}$ we obtain
\begin{equation} \begin{array}{l}\label{Eq_8}  
\sum\limits _{k=1}^{x}\frac{\mu _{k} }{k}  \log \frac{x}{k} +\gamma m_{x} +\frac{1}{2x} M_{x} -\frac{1}{12x^{2} } M_{x} \left(-1\right)+...\approx 1\, ,
\\
m_{x}\log x -\sum\limits _{k=1}^{x}\mu _{k} \frac{\log k}{k} +\gamma m_{x} +\frac{1}{2x} M_{x} -\frac{1}{12x^{2} } M_{x} \left(-1\right)+ ...\approx 1\, ,
\\
\sum\limits _{k=1}^{x}\mu _{k} \frac{\log k}{k}  \approx -1+\left(\log x+\gamma \right) m_{x}+\frac{1}{2x} M_{x} -\frac{1}{12x^{2} } M_{x} \left(-1\right)+ ...\, \, .
\end{array} 
\end{equation} 
Hence it follows that, if $m_{x} =o\left(\frac{1}{\log x} \right)$ at $x\to \infty $ (and, as a consequence, all terms on the right from the term with $m_{x} $ tend to zero), then $\sum _{k=1}^{x}\mu _{k} \frac{\log k}{k}  \to -1$ and vice versa.

For another important case when \textit{s }= 0, $H_{x} \left(0\right)\equiv \left[x\right]$, we have well known formula
\begin{equation}\label{Eq_9}
\sum\limits _{k=1}^{x}\mu _{k} \left[\frac{x}{k} \right] =1,
\end{equation}
with corresponding integral form
\begin{equation}\label{Eq_10}
\int _{1-}^{x}\left[\frac{x}{k} \right] \cdot dM_{y} =1 .
\end{equation}

Using a sieve procedure, the generalized harmonic number can be expanded onto \textit{s-}powers of the consecutive prime numbers $2,\, 3,...,p\le x$ as
\begin{equation}\label{Eq_11}
H_{x} \left(s\right)=1+\frac{1}{2^{s} } H_{\frac{x}{2} }^{\textcircled{\scriptsize 2} } \left(s\right)+\frac{1}{3^{s} } H_{\frac{x}{3} }^{\textcircled{\scriptsize 3} } \left(s\right)+...+\frac{1}{p^{s} } H_{\frac{x}{p} }^{\textcircled{\scriptsize p} } \left(s\right),  
\end{equation}
where by definition recursively
\begin{equation}\label{Eq_12}
H_{x}^{\textcircled{\scriptsize p} } \left(s\right)\equiv H_{x} \left(s\right)-\frac{1}{2^{s} } H_{\frac{x}{2} }^{\textcircled{\scriptsize 2} } \left(s\right)-\frac{1}{3^{s} } H_{\frac{x}{3} }^{\textcircled{\scriptsize 3} } \left(s\right)-...-\frac{1}{p_{-}^{s} } H_{\frac{x}{p_{-} } }^{\textcircled{\scriptsize p\_} } \left(s\right),
\end{equation}
or recurrently
\begin{equation}
\begin{aligned}\label{Eq_13}
H_{x}^{\textcircled{\scriptsize 2} } \left(s\right)&\equiv H_{x} \left(s\right)\, , \\ \\
H_{x}^{\textcircled{\scriptsize 3} } \left(s\right)&\equiv H_{x}^{\textcircled{\scriptsize 2} } \left(s\right)-\frac{1}{2^{s} } H_{\frac{x}{2} }^{\textcircled{\scriptsize 2} } \left(s\right)=H_{x} \left(s\right)-\frac{1}{2^{s} } H_{\frac{x}{2} } \left(s\right)\, , \\ \\
H_{x}^{\textcircled{\scriptsize 5} } \left(s\right)&\equiv H_{x}^{\textcircled{\scriptsize 3} } \left(s\right)-\frac{1}{3^{s} } H_{\frac{x}{3} }^{\textcircled{\scriptsize 3} } \left(s\right)
\\
&=H_{x} \left(s\right)-\frac{1}{2^{s} } H_{\frac{x}{2} } \left(s\right)-\frac{1}{3^{s} } H_{\frac{x}{3} } \left(s\right)+\frac{1}{6^{s} } H_{\frac{x}{6} } \left(s\right)\, , \\
.\,.\,.\,\, ,\\
H_{x}^{\textcircled{\scriptsize p} } \left(s\right)&\equiv H_{x}^{\textcircled{\scriptsize p\_} } \left(s\right)-\frac{1}{p_{-}^{s} } H_{\frac{x}{p_{-} } }^{\textcircled{\scriptsize p\_} } \left(s\right)\, ,
\\
\end{aligned}
\end{equation}
$p_{-} $ is the prime preceding the prime \textit{p}.

Consider asymptotic properties at $x\rightarrow\infty$. Let us apply Euler product formula, which is valid for $Re(s) > 1$ \cite{Bateman04, Montgomery06, Tenenbaum95}, to represent the generalized harmonic number limit as 
\begin{equation}\label{Eq_14}
\begin{aligned}
H_{\infty }\left( s\right)  =\sum_{k=1}^{\infty }\frac{1}{k^{s}}&=\left(
\sum_{a_{2}\geq 0}\frac{1}{2^{a_{2}s}}\right) \cdot \left( \sum_{a_{3}\geq 0}%
\frac{1}{3^{a_{3}s}}\right) \cdot \left( \sum_{a_{5}\geq 0}\frac{1}
{5^{a_{3}s}}\right) \cdot ... \\
&=\prod_{\text{all primes }p}\left( 1-\frac{1}{%
p^{s}}\right) ^{-1}. \\
\end{aligned}
\end{equation}
Further sieving all even number reciprocals yields
\begin{equation}\label{Eq_15}
\begin{aligned}
H_{\infty }^{\textcircled{\scriptsize 3}}\left( s\right)  &=\sum_{\left( k,2\right) =1}^{\infty }%
\frac{1}{k^{s}}=\left( \sum_{a_{3}\geq 0}\frac{1}{3^{a_{3}s}}\right) \cdot
\left( \sum_{a_{5}\geq 0}\frac{1}{5^{a_{3}s}}\right) \cdot ... \\
&=\prod_{\text{all primes }p>2}\left( 1-\frac{1}{p^{s}}\right) ^{-1} \, ,
\\
\frac{H_{\infty }^{\textcircled{\scriptsize 3}}\left( s\right) }{H_{\infty }\left( s\right) }
&=\left( 1-\frac{1}{2^{s}}\right).
\end{aligned}
\end{equation}
Similarly, sieving multiples $p = 3$ we have
\begin{equation}\label{Eq_16}
\begin{aligned}
H_{\infty }^{\textcircled{\scriptsize 5}}\left( s\right)  &=\sum_{\left( k,6\right) =1}^{\infty }%
\frac{1}{k^{s}}=\left( \sum_{a_{5}\geq 0}\frac{1}{5^{a_{3}s}}\right) \cdot
...\cdot \left( \sum_{a_{p}\geq 0}\frac{1}{p^{a_{3}s}}\right) \cdot
...\\
&=\prod_{\text{all primes }p>3}\left( 1-\frac{1}{p^{s}}\right) ^{-1}, \\
\frac{H_{\infty }^{\textcircled{\scriptsize 5}}\left( s\right) }{H_{\infty }\left( s\right) }
&=\left( 1-\frac{1}{2^{s}}\right) \left( 1-\frac{1}{3^{s}}\right).
\end{aligned}
\end{equation}
Continuing the sieving procedure up to any prime $p$ leads to
\begin{equation}\label{Eq_17}
\begin{aligned}
H_{\infty }^{\textcircled{\scriptsize p}}\left( s\right)  &=\sum_{\left( k,p_{-}\#\right) =1}^{\infty
}\frac{1}{k^{s}}=\left( \sum_{a_{p}\geq 0}\frac{1}{p^{a_{3}s}}\right) \cdot
... \\
&=\prod_{\text{all primes }p>p_{-}}\left( 1-\frac{1}{p^{s}}\right) ^{-1}, \\
\frac{H_{\infty }^{\textcircled{\scriptsize p}}\left( s\right) }{H_{\infty }\left( s\right) }
&=\left( 1-\frac{1}{2^{s}}\right) \cdot \left( 1-\frac{1}{3^{s}}\right)
\cdot ...\cdot \left( 1-\frac{1}{p_{-}^{s}}\right).
\end{aligned}
\end{equation}

Consider important case \textit{s }= 1. The sequence of formula transformations for consecutive prime numbers leads to the following set of identities
\begin{equation}\label{Eq_18}
\begin{array}{l} 
H_{x}^{\textcircled{\scriptsize 3} } =H_{x} -\frac{1}{2} H_{\frac{x}{2} } \, =1+\frac{1}{3} +\frac{1}{5} +\frac{1}{7} +\frac{1}{9} +...+\frac{1}{q} \, \, ,\qquad q\le x
\\ \\
{H_{x} -\frac{1}{2} H_{\frac{x}{2} }^{\textcircled{\scriptsize 2} } =1+\frac{1}{3} H_{\frac{x}{3} }^{\textcircled{\scriptsize 3} } +...+\frac{1}{p} H_{\frac{x}{p} }^{\textcircled{\scriptsize p} } =H_{x}^{\textcircled{\scriptsize 3} } ,} 
\\ \\
\frac{H_{x}^{\textcircled{\scriptsize 3} }} {H_{x}} = \left(1-\frac{1}{2} \frac{H_{\frac{x}{2} } }{H_{x} } \right),
\end{array}
\end{equation}
which limit at $x\to \infty $ is
\begin{equation}\label{Eq_19}
\frac{H_{\infty }^{\textcircled{\scriptsize 3} } }{H_{\infty } } =\left(1-\frac{1}{2} \right).
\end{equation}
Similarly for the identities at \textit{p }= 5, we write
\begin{equation}
\begin{array}{l}\label{Eq_20}
{H_{x}^{\textcircled{\scriptsize 5} } =H_{x} -\frac{1}{2} H_{\frac{x}{2} } -\frac{1}{3} H_{\frac{x}{3} } +\frac{1}{6} H_{\frac{x}{6} } =1+\frac{1}{5} +\frac{1}{7} +...+\frac{1}{q} \, , \,\, q\le x} 
\\ \\
\frac{H_{x}^{\textcircled{\scriptsize 5} }} {H_{x}} = \left(1-\frac{1}{2} \frac{H_{\frac{x}{2} } }{H_{x} } -\frac{1}{3} \frac{H_{\frac{x}{3} } }{H_{x} } +\frac{1}{6} \frac{H_{\frac{x}{6} } }{H_{x} } \right)
\end{array}
\end{equation}
with corresponding limit
\begin{equation}\label{Eq_21}
\frac{H_{\infty }^{\textcircled{\scriptsize 5} } }{H_{\infty } } =\left(1-\frac{1}{2} \right)\left(1-\frac{1}{3} \right).
\end{equation}
Continuing the same procedure up to any prime \textit{p}, we have
\begin{equation}\label{Eq_22}
\frac{H_{\infty }^{\textcircled{\scriptsize p}}  }{H_{\infty } } =\left(1-\frac{1}{2} \right)\left(1-\frac{1}{3} \right)...\left(1-\frac{1}{p_{-}} \right).
\end{equation}
Ultimately, sieving all primes \textit{p}, we obtain
\begin{equation}\label{Eq_23}
\frac{1}{\zeta \left(1\right)} =\prod _{{\rm all}\, {\rm primes\; p}}\left(1-\frac{1}{p} \right) =0,
\end{equation}
where
\begin{equation} \label{Eq_24} 
\zeta \left(s\right)=\sum _{n=1}^{\infty }\frac{1}{n^{s}}   
\end{equation} 
is Riemann's zeta function and $\zeta \left(1\right) \equiv H_{\infty }$. Thus we have proved the appropriateness of Euler product formula for $\zeta \left(s=1\right)$.

Consider another important case $s = 0$ and $H_x(0)=[x]$. Using a similar procedure as in \eqref{Eq_11}, $\left[x\right]$ can be represented as
\begin{equation}\label{Eq_25}
\left[x\right]=1+\pi _{\frac{x}{2} }^{\textcircled{\scriptsize 2} } +\pi _{\frac{x}{3} }^{\textcircled{\scriptsize 3} } +...+\pi _{\frac{x}{p} }^{\textcircled{\scriptsize p} },
\end{equation}
where by definition recursively
\begin{equation}\label{Eq_26}
\pi _{x}^{\textcircled{\scriptsize p} } \equiv \left[x\right]-\pi _{\frac{x}{2} }^{\textcircled{\scriptsize 2} } -\pi _{\frac{x}{3} }^{\textcircled{\scriptsize 3} } -...-\pi _{\frac{x}{p} }^{\textcircled{\scriptsize p\_} },
\end{equation}
or recurrently
\begin{equation}\label{Eq_27}
\begin{aligned}
\pi _{x}^{\textcircled{\scriptsize 2} } &\equiv \left[x\right]\, , \\
\pi _{x}^{\textcircled{\scriptsize 3} } &\equiv \pi _{x}^{\textcircled{\scriptsize 2} } -\pi _{\frac{x}{2} }^{\textcircled{\scriptsize 2} } =\left[x\right]-\left[\frac{x}{2} \right]\, , \\
\pi _{x}^{\textcircled{\scriptsize 5} } &\equiv \pi _{x}^{\textcircled{\scriptsize 3} } -\pi _{\frac{x}{3} }^{\textcircled{\scriptsize 3} } =\left[x\right]-\left[\frac{x}{2} \right]-\left[\frac{x}{3} \right]\, +\left[\frac{x}{6} \right]\, \, , \\
{\cdot \cdot \cdot } \\
\pi _{x}^{\textcircled{\scriptsize p} } &\equiv \pi _{x}^{\textcircled{\scriptsize p\_} } -\pi _{\frac{x}{p_{-} } }^{\textcircled{\scriptsize p\_} } \, ,
\end{aligned}
\end{equation}
$p_{-} $ is the prime preceding \textit{p}.

Let us consider the case \textit{p} = 3.
\begin{equation}\label{Eq_28}
\begin{aligned}
\pi _{x}^{\textcircled{\scriptsize 3} }&\equiv \left[x\right]-\left[\frac{x}{2} \right]\, \, =\sum _{(q,2)=1,{\rm \; }q\le x}1 \, \, ,\\
\pi _{x}^{\textcircled{\scriptsize 3} }&=\left[x\right]-\pi _{\frac{x}{2} }^{\textcircled{\scriptsize 2} }=1+\pi _{\frac{x}{3} }^{\textcircled{\scriptsize 3} } +...+\pi _{\frac{x}{p} }^{\textcircled{\scriptsize p} } ,  \\ 
\frac{\pi _{x}^{\textcircled{\scriptsize 3} } }{x}&= \left(\frac{\left[x\right]}{x} -\frac{\left[\frac{x}{2} \right]}{x} \right) ,
\end{aligned}
\end{equation}
which limit at $x\to \infty $ is
\begin{equation}\label{Eq_29}
\mathop{\lim }\limits_{x\to \infty } \frac{\pi _{x}^{\textcircled{\scriptsize 3} } }{x} =\left(1-\frac{1}{2} \right).
\end{equation}
Similarly for the identities at \textit{p }= 5, we write
\begin{equation}\label{Eq_30}
\begin{aligned}
\pi _{x}^{\textcircled{\scriptsize 5} } &\equiv \sum _{(q,6)=1,{\rm \; }q\le x}1 \, \, , \\ 
\frac{\pi _{x}^{\textcircled{\scriptsize 5} } }{x}&=\left(\frac{\left[x\right]}{x} -\frac{\left[\frac{x}{2} \right]}{x} -\frac{\left[\frac{x}{3} \right]}{x} +\frac{\left[\frac{x}{6} \right]}{x} \right)
\end{aligned}
\end{equation}
with corresponding limit
\begin{equation}\label{Eq_31}
\mathop{\lim }\limits_{x\to \infty } \frac{\pi _{x}^{\textcircled{\scriptsize 5} } }{x} =\left(1-\frac{1}{2} \right)\left(1-\frac{1}{3} \right).
\end{equation}
Continuing the same procedure up to any prime \textit{p}, we have
\begin{equation}\label{Eq_32}
\mathop{\lim }\limits_{x\to \infty } \frac{\pi _{x}^{\textcircled{\scriptsize p} } }{x} =\left(1-\frac{1}{2} \right)\left(1-\frac{1}{3} \right)...\left(1-\frac{1}{p\_} \right).
\end{equation}
Thus, from \eqref{Eq_32} and \eqref{Eq_18} we can see that
\begin{equation}\label{Eq_33}
\mathop{\lim }\limits_{x\to \infty } \frac{\pi _{x}^{\textcircled{\scriptsize p} } }{x}= \mathop{\lim }\limits_{x\to \infty } \frac{H_{x}^{\textcircled{\scriptsize p}}  }{H_{x} },
\end{equation}
or

\begin{equation}\label{Eq_34}
\mathop{\lim }\limits_{x\to \infty } \pi _{x}^{\textcircled{\scriptsize p} } \frac{H_{x}}{xH_{x}^{\textcircled{\scriptsize p}}}=1.
\end{equation}

\section{Generalized oscillatory numbers}

The same approach as for generalized number $H_{x} \left(s\right)$ can be applied for corresponding oscillatory number in power \textit{s}
\begin{equation}\label{Eq_35}
M_{x} \left(s\right)=\sum _{k=1}^{x}\frac{\mu _{k} }{k^{s} }.
\end{equation}
Particularly at \textit{s }= 1 and \textit{s }= 0 we have the classic summatory functions \cite{Bateman04,Hardy79,Montgomery06,Tenenbaum95,Tenenbaum01}
\begin{equation}\label{Eq_36}
M_{x} \left(1\right)\equiv m_{x} =\sum _{k=1}^{x}\frac{\mu _{k} }{k} ,
\end{equation}

\begin{equation}\label{Eq_37}
M_{x} \left(0\right)\equiv M_{x} =\sum _{k=1}^{x}\mu _{k} \qquad\text{ -- Mertens' function.}
\end{equation}

Using M\"obius inversion formula \eqref{Eq_2} and \eqref{Eq_3} in the same way as for \eqref{Eq_1}, we have
\[x^{s} \cdot M_{x} \left(s\right)=\sum _{k=1}^{x}\mu _{k} \left(\frac{x}{k} \right) ^{s} ,\] 
\[x^{s} =\sum _{k=1}^{x}\left(\frac{x}{k} \right) ^{s} M_{\frac{x}{k} } \left(s\right). \] 

From where the analog of identity \eqref{Eq_4} is obtained in the form
\begin{equation}\label{Eq_38}
\sum _{k=1}^{x}\frac{1}{k^{s} } \cdot  M_{\frac{x}{k} } \left(s\right)=1,
\end{equation}
also having the integral representation
\begin{equation}\label{Eq_39}
\int _{1- }^{x}M_{\frac{x}{y} } \left(s\right) \cdot dH_{y} \left(s\right)=1.
\end{equation}
Once again, at \textit{s }= 1 and \textit{s }= 0 we have
\begin{equation}\label{Eq_40}
\sum _{k=1}^{x}\frac{1}{k} \cdot  m_{\frac{x}{k} } =1,
\end{equation}
\begin{equation}\label{Eq_41}
\int _{1- }^{x}m_{\frac{x}{y} }  dH_{y} =1
\end{equation}
and
\begin{equation}\label{Eq_42}
\sum _{k=1}^{x}M_{\frac{x}{k} }  =1,
\end{equation}
\begin{equation}\label{Eq_43}
\int _{1- }^{x}M_{\frac{x}{k} }  d\left[y\right]=1,
\end{equation}
respectively.

Obviously, for the equation \eqref{Eq_38}
$$
M_{x} \left(s\right)+\frac{1}{2^{s} } M_{\frac{x}{2} } \left(s\right)+\frac{1}{3^{s} } M_{\frac{x}{3} } \left(s\right)+...+\frac{1}{\left[x\right]^{s} } M_{\frac{x}{x} } \left(s\right)=1
$$
the sieve procedure can be applied. Let us, for example, sieve all \textit{s-}powers of even numbers (sieving \textit{p }= 2). In this case we rewrite \eqref{Eq_38} for $x/2$ and multiply both parts of obtained equation by $\frac{1}{2^{s} } $
\begin{equation}\label{Eq_44}
\frac{1}{2^{s} } M_{\frac{x}{2} } \left(s\right)+\frac{1}{4^{s} } M_{\frac{x}{4} } \left(s\right)+\frac{1}{6^{s} } M_{\frac{x}{6} } \left(s\right)+...+\frac{1}{2^{s} } \frac{1}{\left[\frac{x}{2} \right]^{s} } M_{\frac{\left(x/2\right)}{\left(x/2\right)} } \left(s\right)=\frac{1}{2^{s} } .
\end{equation}
Subtracting \eqref{Eq_44} from \eqref{Eq_38}, we derive 
\begin{equation}\label{Eq_45}
M_{x} \left(s\right)+\frac{1}{3^{s} } M_{\frac{x}{3} } \left(s\right)+\frac{1}{5^{s} } M_{\frac{x}{5} } \left(s\right)+...=1-\frac{1}{2^{s} } .
\end{equation}
Sieving further $p=3$, we have
\begin{equation}\label{Eq_46}
M_{x} \left(s\right)+\frac{1}{5^{s} } M_{\frac{x}{5} } \left(s\right)+\frac{1}{7^{s} } M_{\frac{x}{7} } \left(s\right)+...=\left(1-\frac{1}{2^{s} } \right)\left(1-\frac{1}{3^{s} } \right) .
\end{equation}

Sieve procedure can be continued until any prime \textit{p}, for which $x\ge p{\rm \# }\, \, {\rm (primorial)}\equiv {\rm 2}\cdot {\rm 3}\cdot {\rm 5}\cdot {\rm ...}\cdot p$
\begin{equation}\begin{aligned}\label{Eq_47}
M_{x} \left(s\right)+\frac{1}{p_{+1} ^{s} } M_{\frac{x}{p_{+1} } } \left(s\right)&+\frac{1}{p_{+2} ^{s} } M_{\frac{x}{p_{+2} } } \left(s\right)+...
\\
&=\left(1-\frac{1}{2^{s} } \right)\left(1-\frac{1}{3^{s} } \right)...\left(1-\frac{1}{p^{s} } \right) ,
\end{aligned}
\end{equation}
where $p_{+1} $, $p_{+2} $, etc. are all successive primes after prime \textit{p} up to \textit{x}.

Ultimately, sieving all primes, we obtain
\begin{equation}\label{Eq_48}
M_{\infty } \left(s\right)=\prod _{{\rm all\; primes\; p}}\left(1-\frac{1}{p^{s} } \right) . 
\end{equation}
Hence at \textit{s }= 1 from \eqref{Eq_48} immediately follows the Prime Number Theorem 
\begin{equation}\label{Eq_49}
m_{\infty } \equiv \sum _{k=1}^{\infty }\frac{\mu _{k} }{k}  =\prod _{{\rm all\; primes\; p}}\left(1-\frac{1}{p} \right) =0 .
\end{equation}
At \textit{s }= 0 the sieve procedure for \textit{M(x)} gives
\begin{equation} \begin{array}{l}\label{Eq_50}
{M_{x} +M_{\frac{x}{2} } +M_{\frac{x}{3} } +...+M_{\frac{x}{x} } =1\, ,} 
\\ \\
{M_{x} +M_{\frac{x}{3} } +M_{\frac{x}{5} } +...=0\, ,} 
\\
\bf{. \, . \, . \ \, ,}
\\
{M_{x} +M_{\frac{x}{p_{+1} } } +M_{\frac{x}{p_{+2} } } +...=0\, ,} 
\end{array} 
\end{equation}
for any finite \textit{p} and $x\ge p{\rm \# }$. In this case after sieving all primes, we get uncertainty $0\cdot \infty $ for $M_{\infty }$. 

From Euler product formula, which is valid now at \textit{s} = 1 \eqref{Eq_14}, \eqref{Eq_23} and \eqref{Eq_48}, follows
\begin{equation}\label{Eq_51}
\varsigma \left(s\right)\cdot \vartheta \left(s\right)=1 ,
\end{equation}
where we introduced for symmetry
\begin{equation}\label{Eq_52}
\vartheta \left(s\right)\equiv M_{\infty } \left(s\right)=\sum _{k=1}^{\infty }\frac{\mu _{k} }{k^{s} }.
\end{equation}
For example at $s = 2$ we have
\begin{equation}\label{Eq_53}
M_{\infty } \left(2\right)=\vartheta(2) =\frac{1}{\varsigma \left(2\right)} =\frac{6}{\pi ^{2} }.
\end{equation}
From Riemann's functional equation 
\begin{equation} \label{Eq_54} 
\varsigma \left(s\right)=2^{s} \pi ^{s-1} \sin \left(\pi s/2\right)\Gamma \left(1-s\right)\, \varsigma \left(1-s\right) 
\end{equation} 
and from Euler product follows the functional equations for$\, \, \vartheta \left(s\right)$-function
\begin{equation} \label{Eq_55} 
\vartheta \left(s\right)=\left(2^{s} \pi ^{s-1} \sin \left(\pi s/2\right)\Gamma \left(1-s\right)\right)^{-1} \vartheta \left(1-s\right).
\end{equation} 
\indent In conclusion, using the approach above, we represent some preliminary results for oscillatory part of Chebyshev's \textit{$\psi$} -function

\begin{equation}
\begin{array}{l}\label{Eq_56}
{\log \left[x\right]!\, =\sum\limits _{k=1}^{x}\log k =\sum\limits _{k=1}^{x}\psi _{\frac{x}{k} }  \, \, ,} \\ \\
{\psi \left(x\right)=\sum\limits _{k=1}^{x}\mu _{k} \log  \left[\frac{x}{k} \right]!\, \, \, ,} 
\\ \\
{\sum\limits _{k=1}^{x}\psi _{\frac{x}{k} }  =\sum\limits _{k=1}^{x}\log \frac{k}{x}  +\left[x\right]\log x\, \, ,} 
\\ \\
{\left[x\right]\log x=\sum\limits _{k=1}^{x}\left(\psi _{\frac{x}{k} } +\log \frac{x}{k} \right)\, \, \, , } 
\\ \\
{\psi _{x} +\log x=\sum\limits _{k=1}^{x}\mu _{k}  \left[\frac{x}{k} \right]\log \frac{x}{k} \, \, \, .} 
\end{array} 
\end{equation}
Applying the last equation in \eqref{Eq_56}, we can separate the regular and oscillatory parts of the \textit{$\psi$ }-function
\begin{equation}\label{Eq_57}
\psi _{x} =x\cdot \sum _{k=1}^{x}\mu _{k}  \frac{1}{k} \log \frac{x}{k} -\left(\sum _{k=1}^{x}\mu _{k}  \left\{\frac{x}{k} \right\}\log \frac{x}{k} +\log x\right),
\end{equation}
According to \eqref{Eq_4} and \eqref{Eq_8} the first term in right hand side tends to \textit{x} while the second term is oscillatory part, determined by nontrivial zeros of Riemann's zeta function \eqref{Eq_24} through explicit formula
\begin{equation}\label{Eq_58}
\sum _{\zeta \left(\rho \right)=0}\frac{x^{\rho } }{\rho }  \approx \sum _{k=1}^{x}\mu _{k}  \left\{\frac{x}{k} \right\}\log \frac{x}{k} +\log x .
\end{equation}
More detailed discussions concerning oscillatory parts of the basic functions of prime numbers will be published soon \cite{RAbrarov}.
\bigskip
\bigskip
\bigskip

\noindent---------------------------------------------------------------------\\

\end{document}